\documentclass[12pt]{article}
\usepackage{amsmath}
\usepackage{amssymb}
\usepackage{amsthm}
\usepackage{amscd}
\usepackage{amsfonts}
\usepackage{graphicx}
\usepackage{algorithm}
\usepackage{algorithmic}
\usepackage{setspace}
\newtheorem{theorem}{Theorem}[section]
\newtheorem{lemma}[theorem]{Lemma}

\newtheorem{definition}{Definition}[section]

\newtheorem{remark}{Remark}[section]

\begin{document}
	
\title{\bf On Construction of weighted orthogonal matrices over finite field and its application in cryptography}

\author {Shipra Kumari\thanks {E-mail: shipracuj@gmail.com},  Hrishikesh Mahato\thanks{Corresponding author: E-mail: hrishikesh.mahato@cuj.ac.in }, Sumant Pushp\thanks {E-mail: sumantpushp@gmail.com } \\
	\small Department of Mathematics,
	\small Central University of Jharkhand,
	\small Ranchi-835205, India}
\date{}
\maketitle \setlength{\parskip}{.11 in}
\setlength{\baselineskip}{15pt}
\maketitle
\begin{abstract}
In this article, we propose a method to construct self orthogonal matrix, orthogonal matrix and anti orthogonal matrix over the finite field. Orthogonal matrices has numerous applications in cryptography, so here we demonstrate the application of weighted orthogonal matrix into cryptography. Using the proposed method of construction we see that it is very easy to transmit the private key and can easily convert the encrypted message into original message and at the same time  it will be difficult to get the key matrix for intruder.

\end{abstract}
\textbf{Keywords}: Orthogonal matrix, Self Orthogonal matrix, Anti Orthogonal matrix, Finite field, Encryption scheme \\
\textbf{AMS Subject classification}: 15B10, 11T71, 94A60, 68P25 \\
\hrule

\section{Introduction:}
A square matrix $A$ is said to be Orthogonal matrix of order $n$ if \begin{equation*}
A^TA=AA^T=I_{n}.
\end{equation*}
This leads that the inverse of an Orthogonal matrix is its transpose.\\ 
P.G. Farrell introduced a very interesting concept which is "anti-code" \cite{Farrell}. Generally a "code" is designed to have a large minimum distance between its code words. As opposite of "large minimum distance" is "small maximum distance". Farrell used the term "anti-code" to describe a set of $n$-tuple designed to have small maximum distance between its code words \cite{selforthogonal}. To get inspired by Farrel's "anti code" James L. Massey \cite{selforthogonal} introduced anti-orthogonal matrix, as opposite of $I$ is $-I$. If $BB^T= -I$  then $B$ is known as anti-orthogonal matrix. If  each row of a matrix is orthogonal to every row including itself  is self orthogonal matrix  i.e. if $CC^T= \textbf{0}$ , then $C$ is a self-orthogonal matrix.\\
\par Cryptography is the study of techniques through which a secured communication of message or information executed i.e. communicated message/information can’t be understood by anyone who is not the intended receiver. The method to obtain the cipher text from plain text is encryption and the reverse method is decryption. There are several techniques for encryption scheme \cite{encrypt1, encrypt2, hill, koukouencrypt}. In $1929$ Lester S. Hill \cite{hill} introduced the Hill cipher in which  an invertible matrix is  used as a private key and inverse of that matrix is used to decrypt. It is very difficult to obtain the inverse of a matrix for large size. Later on orthogonal matrices are used as a private key to disposed off the calculation of inverse of the matrix. As inverse of an orthogonal matrix is its transpose. C. Koukouvinos and D.E. Simos \cite{koukouencrypt} also developed an encryption scheme using circulant
Hadamard core in which the first row of the Hadamard core required to be transmitted as a
private key. 
 \par In this paper we present a method to construct the  weighted orthogonal matrices over a finite field $GF(q)$, where $q=p^{\alpha}$ such that $p$ is a prime and $\alpha$ is a positive integer. Also we develop an encryption scheme using the  weighted orthogonal matrices over finite field $GF(q)$. Our main aim is to protect the information about the elements of the key matrix. In this encryption scheme there is no need to transmit the matrix or first row  or any element of the matrix to get the private key. Along with the encrypted message we need to transmit three numbers which include the size of the key matrix and additionally a primitive polynomial in case of $ \alpha >1$. So for intruder it will be difficult to find the private key and the numbers make an easy transmission.  \\
The paper is organized as follows:\\
\par In section $2$, the required definitions and information have been discussed as preliminaries.    In section $3$, we present the methods to construct orthogonal, self-orthogonal and anti-orthogonal matrices of order $q$ over the finite field $GF(q)$. In section $4$, we present the application of orthogonal matrix into cryptography. Here encryption scheme and analysis of time complexity of the scheme have been discussed. In section $5$, we provide the example of the encryption scheme defined in section $4$ and in section $6$ we conclude the results.
  
\section{Preliminaries}
\begin{definition}
	\textbf{Weighted orthogonal matrix over a field}\\
	 A square matrix $A=[a_{ij}]_{n\times n}$, where $a_{ij}\in GF(q)$ is said to be an orthogonal matrix with weight $k$ if \begin{equation}\label{def:orthogonality}
	 AA^T=kI_{n}
	 \end{equation}  over $GF(q)$.\\
	 	In particular if $k=0,\ 1, \ -1 $, then  A is a self orthogonal matrix, orthogonal matrix and anti orthogonal matrix of order $n$ respectively over $GF(q)$
	  \begin{remark}
	  	If $q=3$, then the matrix defined over $GF(3)$ which satisfies the equation \eqref{def:orthogonality} is a weighing matrix of order $3$.
	  \end{remark}	 
\end{definition} 

	
\begin{definition} \textbf{Primitive root}\\
Let $GF(q)$ be a finite field of order $q$. A primitive root is an element of $GF(q)$ which is a  generator of the multiplicative group of the field $GF(q)$.

\end{definition}
\section{Construction of Orthogonal Matrix}
In this section we discuss some method to construct the orthogonal matrices of order $q$ over $GF(q)$.\\
Let $q=p^\alpha$, where $p$ is a prime and $\alpha$ is any positive integer,
\begin{equation*}
GF(q) = \left\{0= a_{0}, \ a_{1}, \ a_{2},  \cdots, a_{q-1} \right\}
\end{equation*}
\begin{lemma}\label{lemma:divisibility}
	Let $q=p^{\alpha}$, where $p$ is a prime and $\alpha$ is any positive integer. Then for any $1\leq k<q-1$  \begin{equation}
	\sum_{i=0}^{q-1}a_{i}^k \equiv 0(\bmod \ p).	\end{equation}
\end{lemma}
\textit{Proof.}	Let $g$ be a primitive root of $GF(q)$ so  for every non-zero element
	$a_{i} \in GF(q)$, $a_{i}$ can be written as 
	\begin{equation*} a_{i}=g^\beta \ ; \ \ \ \ \ \ \ \ \ \ \ \ \  1 \leq \beta \leq q-1
	\end{equation*}
	Now \begin{align*}
	 \sum_{i=0}^{q-1}a_{i}^k &= 0+\sum_{\beta=1}^{q-1}(g^\beta)^k \\
	 &=g^k+g^{2k}+\cdots+g^{(q-1)k}\\
	& = g^k\{ 1+g^k + g^{2k}+ \cdots + g^{(q-2)k}\}\\
	&=g^k \frac{\{ (g^k)^{(q-1)}-1\}}{g^k-1}\\
	&=g^k \frac{\{(g^{(q-1)})^k-1\}}{g^k-1}\\
	\end{align*}
	 As $g$ is a primitive root of $GF(q)$ and $1 \leq k < q-1$  so $g^{(q-1)} \equiv 1(\bmod \ p)$  and $g^k\neq 1$.\\
	 Thus \begin{equation*}
	 \sum_{i=0}^{q-1}a_{i}^k \equiv 0(\bmod \ p)
	 \end{equation*}
	 \qed	

\begin{theorem}\label{Theorem:AA^T=0}Let $q=p^{\alpha}$, where $p$ is a prime and $A$ is a square matrix of order $q$ with entries $[a_{ij}] \equiv (a_{j}^{t}-a_{i}^{t})(mod \ p)$ where $0 \leq i,j \leq q-1 ; 1 \leq t < \frac{q-1}{2}$ and $a_{i}, a_{j} \in GF(q)$.
 Then  $A$ is a self orthogonal skew symmetric matrix over $GF(q)$.
\end{theorem}
\textit{Proof.}
Let  $A$ be a square matrix of order $q$ with entries $[a_{ij}] \equiv (a_{j}^{t}-a_{i}^{t})(mod \ p)$,  where $0 \leq i,j \leq q-1$ and $1 \leq t < \frac{q-1}{2}$.\\
So, each element of $i^{th}$ row can be expressed as 
\begin{equation*}
a_{ij}= (a_{j}^t-a_{i}^t)(\bmod \ p) \hspace{2cm} \mbox{for $0 \leq j \leq q-1$}
\end{equation*}
Now,  inner product of $u^{th}$ row with itself is
	 \begin{align*}
	 \left<R_{u} \cdot R_{u}\right> & \equiv \sum_{j=0}^{q-1} (a_{j}^t-a_{u}^t)\cdot (a_{j}^t-a_{u}^t) (\bmod \ p)\\
	 & \equiv \sum_{j=0}^{q-1} \left\{a_{j}^{2t} -a_{j} ^ta_{u}^t-a_{u}^ta_{j}^t+a_{u}^ta_{u}^t \right\}(\bmod \ p)\\
	 & \equiv \left\{\sum_{j=0}^{q-1}a_{j}^{2t}- a_{u}^t\sum_{j=0}^{q-1}a_{j}^t-a_{u}^t\sum_{j=0}^{q-1}a_{j}^t+\sum_{j=0}^{q-1}a_{u}^{2t} \right\}(\bmod \ p)\\
	 & \equiv \left\{\sum_{j=0}^{q-1}a_{j}^{2t}- a_{u}^t\sum_{j=0}^{q-1}a_{j}^t-a_{u}^t\sum_{j=0}^{q-1}a_{j}^t+qa_{u}^{2t}\right\}(\bmod \ p)\\
	  \end{align*}
Since $1 \leq t < \frac{q-1}{2} \Rightarrow 2t < q-1$, so using lemma \eqref{lemma:divisibility}, each sum of right hand side is zero. So $\left<R_{u} \cdot R_{u}\right>\equiv \textbf{0} \ (\bmod \ p)$.\\

 Now, inner product of $u^{th}$ row $v^{th}$ row of the matrix $A$\\
\begin{align*}
\left<R_{u} \cdot R_{v} \right> & \equiv \sum_{j=0}^{q-1}(a_{j}^t-a_{u}^t)\cdot(a_{j}^t-a_{v}^t)(\bmod \ p)\\
& \equiv \sum_{j=0}^{q-1}\left\{ a_{j}^{2t} - a_{j}^ta_{v}^t- a_{u}^ta_{j}^t+a_{u}^ta_{v}^t \right\}(\bmod \ p)\\
&\equiv \left\{\sum_{j=0}^{q-1} a_{j}^{2t}-a_{v}^t \sum_{j=0}^{q-1}a_{j}^t - a_{u}^t\sum_{j=0}^{q-1}a_{j}^t+\sum_{j=0}^{q-1}a_{u}^ta_{v}^t \right\}(\bmod \ p)\\
&\equiv \left\{\sum_{j=0}^{q-1} a_{j}^{2t}-a_{v}^t \sum_{j=0}^{q-1}a_{j}^t - a_{u}^t\sum_{j=0}^{q-1}a_{j}^t+q\cdot a_{u}^ta_{v}^t\right\}(\bmod \ p)
\end{align*}
Since $1 \leq t < \frac{q-1}{2} \Rightarrow 2t < q-1$, so using lemma \eqref{lemma:divisibility}
\begin{align*}
\sum_{j=0}^{q-1} a_{j}^{2t} \equiv(\bmod \ p) \hspace{1cm} \mbox{and} \hspace{1cm}
\sum_{j=0}^{q-1} a_{j}^{t} \equiv 0(\bmod \ p)
\end{align*}
So, 
\begin{align*}
\left<R_{u}\cdot R_{v} \right> & \equiv 0(\bmod \ p)
\end{align*}
So, we have $AA^{T}= \textbf{0} \ (\bmod \ p)$ i.e. $A$ is a self orthogonal matrix of order $q$.\\
Now, 
to show that $A$ is a skew symmetric matrix it is enough to show that $a_{ij}=-a_{ji}$.
\begin{align*}
a_{ij}& \equiv (a_{j}^{t}-a_{i}^{t})(\bmod \ p)\\
&\equiv -(a_{i}^{t}-a_{j}^{t})(\bmod \ p)\\
&=-a_{ji}
\end{align*} 
So, $A$ is a skew symmetric matrix of order $q$.
\qed

In particular, if matrix $A$ has same sequence of leading rows and columns then it can be observed that all diagonal entries of $A$ are $0$.

\begin{theorem}\label{theorem:A+rI}
	If  $B$ is a self orthogonal skew symmetric  matrix over $GF(q)$ with diagonal $0$ then $B+rI$, where $ r  \in GF(q)$, is an orthogonal matrix with weight $r^2$, a quadratic residue in $GF(q)$.
\end{theorem}
\textit{Proof.}
	Here
	\begin{align*}
	(B+rI)(B+rI)^T &=(B+rI)(B^T+rI)\\
	&= BB^T+rB+rB^T+r^2I	
	\end{align*}
Since $B$ is a self orthogonal skew symmetric matrix so $BB^T= \textbf{0}$ and $B^T=-B$ \\
Hence, \begin{equation*}
(B+rI)(B+rI)^T=r^2I
\end{equation*}  So $B+rI$ is an orthogonal matrix with weight $r^2$.
\qed
\begin{remark}
For any  $q=p^{\alpha}$, $1$ is a quadratic residue in $GF(q)$. So  $B+rI$ is an orthogonal matrix of order $q$ for $r^2 \equiv 1(\bmod \  p)$. However, for $q=p^{\alpha}$ with $p \equiv 1(\bmod \ 4)$, $-1$ is a quadratic residue in $GF(q)$ implies that $B+rI$ is an anti-orthogonal matrix of order $q$ for $r^2 \equiv -1 (\bmod \ p)$.

\end{remark}

From the above discussion we observed that the weight of the orthogonal matrices is a quadratic residue of $GF(q)$. It is know that each element of $GF(2^\alpha)$ is a quadratic residue in $GF(2^\alpha)$ \cite{sum of two squares}. So for each $k \in GF(2^\alpha)$ there exist an orthogonal matrix of order $2^\alpha$ with weight $k$.

We are interested to find orthogonal matrices of any weight $k \in GF(q)$ over $GF(q)$, where $q$ is an odd prime power.\\
It is  known that every element of a finite field $GF(q)$ can be written as sum of two squares \cite{sum of two squares}. Using this result the orthogonal matrices of order $2q$ of weight $k \in GF(q)$ may be obtained in next results.

 \begin{theorem}\label{theorem:AB}
 	If $A$ and $B$ are two self orthogonal matrices of order $q$ obtained by theorem \eqref{Theorem:AA^T=0} then $AB \equiv 0(\bmod \ p)$.	
 \end{theorem}
\textit{Proof.}
	Let $A=[a_{ij}]$ and $B=[b_{ij}]$ be two self orthogonal matrices over $GF(q)$. Then using theorem \eqref{Theorem:AA^T=0} 
	\begin{equation*}
	a_{ij} = (a_{j}^t-a_{i}^t)(\bmod \ p)
	\end{equation*}
	and
	\begin{equation*}
	b_{ij} = (a_{j}^s-a_{i}^s)(\bmod \ p)
	\end{equation*} 
	where $ 1 \leq \ t, \ s < \frac{q-1}{2}$ and $ a_{i} \in GF(q)$ for $ 0 \leq i, \ j \leq q-1$\\
	Now, 
	\begin{align*}
	AB & = \left\{\sum_{j=0}^{q-1} a_{ij}b_{jk}\right\}(\bmod \ p)\\
	& = \left\{\sum_{j=0}^{q-1} \left(a_{j} ^t-a_{i}^t\right) \left(a_{k}^s-a_{j}^s\right)\right\} (\bmod \ p)\\
	& = \left\{\sum_{j=0}^{q-1} a_{j}^t \cdot a_{k }^s -a_{j} ^t \cdot a_{j}^s- a_{i} ^t \cdot a_{k} ^s+ a_{i}^t \cdot a_{j}^s \right\} (\bmod \ p)\\
	&= \left \{a_{k}^s\sum_{j=0}^{q-1}a_{j}^t -\sum_{j=0}^{q-1}a_{j}^{t+s} - \sum_{j=0}^{q-1}a_{i}^t \cdot a_{k}^s + a_{i}^t \sum_{j=0}^{q-1}a_{j}^s \right\} (\bmod \ p)\\
	& =\left\{a_{k}^s\sum_{j=0}^{q-1}a_{j}^t -\sum_{j=0}^{q-1}a_{j}^{t+s} - q \cdot a_{i}^t \cdot a_{k}^s + a_{i}^t \sum_{j=0}^{q-1}a_{j}^s \right\} (\bmod \ p)\\
	\end{align*}

	As, $ 1 \ \leq \ t, \ s \ < \frac{q-1}{2}$, \ \  so \  $t+s < q-1$. Now using lemma \eqref{lemma:divisibility}
	\begin{align*}
	\sum_{j=0}^{q-1}a_{j}^t \equiv 0(\bmod \ p) \ ,  \ \sum_{j=0}^{q-1}a_{j}^s \equiv 0(\bmod \ p) \ \  \mbox{and} \ \sum_{j=0}^{q-1}a_{j}^{t+s} \equiv 0(\bmod \ p)
	\end{align*}
	Hence \begin{equation*}
	AB \equiv 0(\bmod \ p)
	\end{equation*}
	\qed

\begin{theorem}
	If $A$ and $B$ are two skew symmetric self orthogonal matrices of order $q=p^{\alpha}$, where $p$ is a prime, then there is a weighted orthogonal matrix of order $2q$ with any weight $k \in GF(q)$.
\end{theorem}
\textit{Proof.}
	Since $A$ and $B$ are skew symmetric self orthogonal matrices of order $q$. Since each element of $GF(q)$ can be expressed as a sum of two squares, consider $k=r^2+s^2$ for some $r, \ s \in GF(q)$.\\
	Consider
	
	\begin{equation}
	A^{'}= \left[\begin{array}{cc}
	A+rI & B+sI\\
	-(B+sI)^{T} & (A+rI)^{T}
	\end{array}
	\right]
	\end{equation}

Then
	\begin{align*}
	A'A'^{T}&=\left[\begin{array}{cc}
	A+rI & B+sI\\
	-(B+sI)^{T} & (A+rI)^{T}
	\end{array}
	\right]\left[\begin{array}{cc}
	A+rI & B+sI\\
	-(B+sI)^{T} & (A+rI)^{T}
	\end{array}
	\right]^{T}\\
	&=\left[\begin{array}{cc}
	A+rI & B+sI\\
	-(B+sI)^{T} & (A+rI)^{T}
	\end{array}
	\right]\left[ \begin{array}{cc}
	A^T+rI & -B-sI\\
	B^T+sI & A+rI
	\end{array}
	\right]\\
&=	\left[ \begin{array}{cc}(r^2+s^2)I & -AB+BA\\
	-B^TA^T+A^TB^T & (r^2+s^2)I
	\end{array}
	\right]\\ 
	&=\left[\begin{array}{cc}
	(r^2+s^2)I & \textbf{0}\\ 
	\textbf{0} & (r^2+s^2)I
	\end{array}
	\right]  \hspace{3cm} \mbox{since $AB=BA=0(\bmod \ p)$}\\
	&= (r^2+s^2)I\otimes I_{2}\\
	&=kI_{2q}
	\end{align*}
	Therefore $A'$ is a weighted orthogonal matrix of order $2q$ with weight $k$.
	\qed

In particular, if $r=1$ and $s=0$, \  $A'$ is an orthogonal matrix of order $2q$.
\begin{theorem}
	Let $A+rI$ be a weighted orthogonal matrix of order $p$ with weight $r^2$. If $H_{4n}$ is a Hadamard matrix of order $4n$ then the kronecker product of $(A+rI)$ and $H_{4n}$ is a weighted orthogonal matrix with weight $4nr^2$.
\end{theorem}
\textit{Proof.}
Let $A+rI$ be a weighted orthogonal matrix of order $p$ with weight $r^2$. So,
\begin{equation*}
(A+rI)(A+rI)^T= r^2I
\end{equation*}
Since $H_{n}$ is a Hadamard matrix of order $4n$ so by the definition of Hadamard matrix
\begin{equation*}
H_{4n}H_{4n}^T= 4nI_{4n}
\end{equation*}
Now 
\begin{align*}
&\left \{ H_{4n} \otimes (A+rI)\right\}\left\{ H_{4n} \otimes (A+rI)\right\}^T  \\
&=\left \{ H_{4n} \otimes (A+rI)\right\}\left \{ H_{4n}^T \otimes (A+rI)^T\right\}\\
& = H_{4n}H_{4n}^T \otimes (A+rI)(A+rI)^T\\
&= 4nI_{4n} \otimes r^2I_{p}\\
&=4nr^2I_{4np}
\end{align*}
This shows that kronecker product of  Hadamard matrix and weighted orthogonal matrix is again a weighted orthogonal matrix.
\qed

 \textbf{Note:} If $ \gcd(n, p) \neq 1$ then $H_{4n}\otimes (A+rI)$ is a self orthogonal matrix of order $4np$.
 
\section{Encryption Scheme}

Orthogonal matrix is widely used in cryptography. In this section we have used the weighted orthogonal matrix, discussed  above, in encryption scheme. \\
Suppose the language, in which message or information is written, contains $q=p^\alpha$ ($p$ is a prime) distinct characters. Let $M$ be the numeric plain text of the corresponding message which has to be transmitted securely. Here weighted orthogonal matrix of order $q$ is used to encrypt the message. 
 The encrypted message $C$ is obtained from the plain text $M$ using the transformation
\begin{equation}
C \equiv WM (\bmod \ p)
\end{equation}
where $W$ is a  weighted orthogonal matrix of order $q$ with weight $r^2$. $W$ is defined as 
\begin{equation}
W= A+rI
\end{equation}
where $A$ is a self orthogonal matrix of order $q$. \\
Using theorems~\eqref{Theorem:AA^T=0}~and~\eqref{theorem:A+rI} the key matrix $W$ of order $q=p^\alpha$ over a field $GF(q)$ may be constructed. For $\alpha=1$ the field may be considered as $\mathbb{Z}_q$ but for $\alpha >1$ it requires a primitive polynomial $P(p^\alpha)$ to construct the required field $GF(q)$. Thus formation of private key for $\alpha=1$ it requires the numbers $(p, \ t, \ r)$ however for $\alpha >1$ it requires $P(p^\alpha)$ as an additional information.\\
Here algorithm has been discussed to develop the key matrix

\begin{algorithm}		
	\caption{Formation of key matrix } \label{algorithm:key}
	\begin{algorithmic}[1]
		\REQUIRE $p, \ \alpha, \ t, \ r, \ P(p^{\alpha})$
		\ENSURE $1 \leq t < \frac{q-1}{2}$
		\IF {$\alpha =1$}
		\STATE $q \leftarrow p$
		\STATE Form an array $\mathbb{Z}_{q}$ and\\
		$\mathbb{Z}_{q} = [0, \ 1, \ 2, \, \cdots, \ \ q-1]$ \hspace{3cm} element of the field $GF(q)$
		\FOR{$i \gets 0$ to $q-1$ }
		\STATE $ Z \leftarrow \left\{ a_{i} ^t \ (\bmod p) : a_{i}\in \mathbb{Z}_{q} \ \ \ \mbox{and} \ \ \ 0\leq i \leq q-1 \right\}$
		
		\ENDFOR

		\ELSE[$\alpha >1$]
		\STATE $ q \leftarrow p^{\alpha}$
		\STATE Form an array $F$
		and element of $F$ is defined as:\\
		$F = \left\{\sum_{i=0}^{\alpha -1}(b_{i}x^{i} ) : b_{j} \in \mathbb{Z}_{p}\right\}$\\
		
		\FOR {$i \gets 0 $ to $q-1$}
		\STATE $Z \leftarrow \left\{ a_i^t \left(\bmod P(p^{\alpha})\right): a_i \in F
		\right\}$
		\ENDFOR
		\ENDIF
		\STATE Construct a square matrix $A=[e_{ij}]$ of order $q$
		\FOR{$ i \gets 0$ to $q-1$}
		\FOR {$j \gets 0$ to $q-1$}
		\STATE $ e_{ij} \leftarrow (a_{j}^t-a_{i}^t)(\bmod \ p): a_j^t,\ a_i^t \in Z $
		\ENDFOR
		\ENDFOR
	\end{algorithmic}
\end{algorithm}
 The encryption algorithm for the scheme is given by:

\begin{algorithm}
	\caption{Encryption Algorithm}
	\begin{algorithmic}[1]
		\REQUIRE Require msg to encrypt\\
		\STATE select $(p, \ \alpha,\ t, \ r, P(p^\alpha))$\\
		\IF {$\alpha=1$}
		\STATE$k \leftarrow (p, \ t, \ r)	 \hspace{7.8cm}$  set private key\\
		\ELSE[$\alpha >1$]
		\STATE$k \leftarrow (p, \ t, \ r, P(p^\alpha))	 \hspace{6.8cm}$  set private key\\
		\ENDIF
		\STATE $W \leftarrow A+rI$ \hspace{6.2cm} construct matrix $A$ using algorithm \eqref{algorithm:key}
		\STATE $M \leftarrow$ convert(msg) \hspace{5cm} convert msg into its numerical value\\
		\STATE $C \leftarrow WM (\bmod p)$
		\STATE Transmit$(k,C)$
	\end{algorithmic}
\end{algorithm}

In general, to obtain the original message it requires $W^{-1}$. Since $W$ is a weighted orthogonal matrix, so $W^{-1}= lW^T$, where $l$ is a solution of $ r^2x \equiv 1(\bmod \ p)$. The intended receiver has to decrypt the message using the transformation
\begin{equation}\label{transform:M}
M \equiv lW^TC (\bmod \ p)
\end{equation}

Along with the security, the main objective of the cryptography is to decrypt the message uniquely.
\begin{theorem}
	If $C$ is an encrypted message with the encryption algorithm then the decrypted message  
	\begin{equation*}
	D \equiv lW^TC(\bmod \ p)
	\end{equation*}is uniquely determined and is equal to $M$, where $l$ is a solution of \ $r^2 x \equiv 1(\bmod \ p)$
\end{theorem}
	Let $C$ be an encrypted message with transformation 
	\begin{equation*}
	C\equiv WM(\bmod  p)
	\end{equation*}
	Since 
	\begin{align*}
	D &\equiv lW^TC(\bmod p) \\
	&\equiv lW^TWM (\bmod p)\\
	&\equiv lr^2IM(\bmod p) \hspace{1cm}\mbox{as $W$ is a weighted orthogonal matrix so $W^TW=r^2I$}\\
	& \equiv lr^2M(\bmod p)\\
	& \equiv M(\bmod p) \hspace{2cm} \mbox{as $l$ is a solution of $ r^2 x \equiv 1 (\bmod p)$}
	\end{align*}
	It shows that message is uniquely decrypted.
The decryption algorithm is defined below
\begin{algorithm}
	\caption{Decryption Algorithm}
	\begin{algorithmic}[1]
		\REQUIRE Require obtained cipher text\\
		\STATE receive $(C, \ p, \ t, \ r, P(p^\alpha))	$\\
		\STATE obtain the value of $l$ \hspace{3cm} \mbox{*// using euclidean algorithm}
		\STATE construct matrix $A$ \hspace{3.3cm} \mbox{*// using algorithm \eqref{algorithm:key}}
		\STATE $W \leftarrow A+rI$ 
		\STATE $M \leftarrow lW^TC (\bmod \ p)$ \hspace{4cm} \mbox{*// decrypt the message}
	\end{algorithmic}
\end{algorithm}

\subsection{Cryptanalysis of Encryption scheme}

In the above discussed encryption scheme, sender transmits  numbers $(p,  \ r, \ t)$ and $P(p^\alpha)$ in either case as an additional information to form the private key. To get the plain text intended receiver has to use the  transformation \eqref{transform:M}.
So receiver has to find out $W^T$ and $l$ only. The matrix $W^T$ and integer $l$  may be obtained by using algorithm \eqref{algorithm:key} and the Euclidean algorithm respectively. The time complexity to find $W^T$ and the integer $l$ is $\mathcal{O}(q^3)$ and $\mathcal{O}(log \ q)$ the respectively.\\
The main aim of encryption scheme is to protect the information about the private key. But the intruder always tries all possible combinations to find the keys and checks which one of them return the plain text . Since in this encryption scheme numbers $(p, \ r, \ t)$ and $P(p^\alpha)$ in either case as an additional information  will be transmitted as a private key. So it will be difficult to guess the size of the key matrix. And if any one could guess the order of a key matrix  it is still difficult to obtain the key matrix $W$. As matrix $W$ consist the elements of $GF(q)$. So the size of the key space, $K(W)$, is $q^{q^2}$. Thus the computational complexity to find the key matrix $W$ is $\mathcal{O}(q^{q^2})$ and it increases exponentially.\\

\textbf{Observation}
\begin{enumerate}
	\item To form the private key numbers $(p, \ r, \ t, P(p^\alpha))$ is used, so it is easy to transmit.
	\item As orthogonal matrix is used to encrypt the message so it is easy to decrypt the message as inverse of an orthogonal matrix is its transpose.
	\item For intruder after knowing the actual size still it is difficult to find the key matrix as  the computational complexity to find the key matrix is $\mathcal{O}(q^{q^2})$.
\end{enumerate}
\section{Example}
Suppose the plain text  COVID$-19$ has to be transmit. Firstly we convert the plain text into its corresponding numerical value (ASCII Code)
\begin{equation*}
\begin{array}{|c|c|c|c|c|c|c|c|}
\hline
	C & O & V & I & D & - & 1 & 9\\
	\hline
	67 & 79 & 86 & 73 & 68 & 45 & 49 & 57\\
	\hline
	\end{array}
	\end{equation*}
	Here highest value present in the plain text is $86$. Consider $p=89$ (we need to choose a prime number which is greater or equal to the highest numeric value present in plain text). In plain text we need to add "space" $81$ times. Thus 
	\begin{equation*}
	M= \left[ \begin{array}{ccccccccccccc}
	67 & 79 & 86 & 73 & 68 & 45 & 49 & 57 & 32  \underbrace{\cdots}_{81 \ times} & 32 & 32\\
	\end{array}
	\right]^T_{1 \times 89}
	\end{equation*}
	Consider $r=5, t=2$ and construct the weighted orthogonal matrix $W$ of order $89$ by applying  algorithm \eqref{algorithm:key}.
	The cipher text $C$ is obtained by the transformation 
	\begin{equation*}
	C \equiv WM(\bmod \ 89)
	\end{equation*}
	is 
	\begin{equation*}
	C= \left[
	\begin{array}{ccccccccccccccccccccccccccccccccccccccccccccccccccccccccccccccccccccccccccccccccccccccccccccccc}
56 & 26 & 58 & 77 & 45 & 10 & 19 & 46 & 84 & 67 & 48 & 27 & 4 & 68 & 41 & 12 & 70 & 37 & 2 & 54 & 15 & 63\\ & 20 & 64 & 17 & 57 & 6 & 42 & 76 & 19 & 49 & 77 & 14 & 38 & 60 & 80 & 9 & 25 & 39 & 51 & 61 & 69 & 75\\ & & 79 & 81 & 81 & 79 & 75 & 69 & 61 & 51 & 39 & 25 & 9 & 80 & 60 & 38 & 14 & 77 & 49 & 19 & 76 & 42 \\ && & 6  & 57 & 17 & 64 & 20 & 63 & 15 & 54 & 2 & 37 & 70 & 12 & 41 & 68 & 4 & 27 & 48 & 67 & 84\\ &&& & 10 & 23 & 34 & 43 & 50 & 55 & 58 
	\end{array}
		\right]^T_{1 \times 89}
	\end{equation*}
	
	It is known that $ 0-31$ and $127 $ are non printable character. So for our convenience we define the corresponding  numerical value (ASCII code) 
	\begin{align*}
	0 & \rightarrow (0)*\\
	1& \rightarrow (1)*\\
	2 &\rightarrow (2)*\\
	& \vdots	
	\end{align*}
	Convert the numeric value of $C$ in their corresponding character
	\begin{equation*}
	C= \left[ 
	\begin{array}{ccccccccccccccccccccccccccccccccccccccccccccccccccccccccccccccccccccccccccccccccccccccccccccccc}
	8 & (26)* & : & M & - & (10)* & (19)* & . & T & C & 0 & (27)* & (4)* & D & )\\ & (12)* & F & \% & (2)* & 6 & (15)* & ? & (20)* & @ & (17)* & 9 & (6)* & * & L\\ & (19)* & 1 & M & (14)* & \& & < & P & (9)* & (25)* & ' & 3 & = & E &K \\ & O & Q & Q &O & K & E & = & 3 & ' & (25)*  & (9)* & P & < & \& \\ & (14)* & M & 1 & (19)* & L & * & (6)* & 9 & (17)* & @ & (20)*  & ? & (15)* & 6\\ & (2)* & \% & F & (12)* & ) & D & (4)* & (27)* & 0 & C & T & (10)* & (23)* & " \\ & + & 2 & 7 & : 
	\end{array}
	\right]^T_{1 \times 89}
	\end{equation*}
	We transmit $C, 89, 5, 2$ .\\
	To obtain the plain text intended receiver has to construct the weighted orthogonal matrix $W$ using algorithm \eqref{algorithm:key}.
	Then by the help of transformation 
	\begin{equation*}
	M= lW^TC(\bmod \ p)
	\end{equation*} 
	$l$ is the solution of $25^2x \equiv 1(\bmod \ 89)$ as $r=5$ so in this case $l=57$ and the plain text will be
		\begin{equation*}
	M= \left[ \begin{array}{ccccccccccccc}
	67 & 79 & 86 & 73 & 68 & 45 & 49 & 57 & 32  \underbrace{\cdots}_{81 \ times} & 32 & 32\\
	\end{array}
	\right]^T_{1 \times 89}
	\end{equation*}

	i.e. 
	\begin{equation*}
	\begin{array}{|c|c|c|c|c|c|c|c|}
	\hline
	C & O & V & I & D & - & 1 & 9\\
	\hline
	67 & 79 & 86 & 73 & 68 & 45 & 49 & 57\\
	\hline
	\end{array}
	\end{equation*}

\section{Conclusion}
In this paper we have developed a method to construct weighted orthogonal matrices over a finite field $GF(q)$. There is no non-zero self-orthogonal matrix over $\mathbb{R}$ however self orthogonal matrix has been developed over $GF(q)$ in modular arithmetic. Furthermore an encryption scheme has been developed using the obtained weighted orthogonal matrices. In this scheme it is easier to transmit a finite number of numbers which yields private key.


\begin{thebibliography}{99}
	\bibliographystyle{}
	
	\bibitem{koukouencrypt}C. Koukouvinos, and D E. Simos, \textit{Encryption schemes based on Hadamard Matrices with Circulant Cores}, Journal of Applied Mathematics and Bioinformatics, Vol. $3$, pp $17-41$, $(2013)$.
\bibitem{anticode}	F. J. MacWilliams and N. J. A. Sloane, \textit{ The Theory of Error Correcting Codes} Amsterdam: North-Holland,  pp $548-556$, $(1977)$.
	
\bibitem{singer}J. Singer, \textit{ A theorem in finite projective geometry and some applications
	to number theory}, Trans. Amer. Math. Soc., Vol. $43$, pp $377-385$, $(1938)$.
	
\bibitem{selforthogonal} J. L. Massey, \textit{Orthogonal, Anti orthogonal and Self-Orthogonal
Matrices and their Codes}, Signal \& Information Processing Laboratory
Swiss Federal Institute of Technology
CH-8092 Zurich, Switzerland	
	
	\bibitem{} J. M. Goethals, and J. J. Seidel, \textit{Orthogonal Matrices with  Zero Diagonal}, Geometrics and Combinatorics, pp $257-266$, $(1991)$,
		
	
	
\bibitem{hill}	Lester S. Hill, \textit{ Cryptography in an Algebraic Alphabet}, The American Mathematical Monthly,
	Vol. $36$, No. $6$, $306- 312$ (Jun-July, $1929$).
	
	\bibitem{mar} Marshall, Jr. Hall, \textit{Combinatorial Theory}, Wiley Interscience Series in Discrete Mathematics, MAc-Wiley $(1986)$.
	
	
\bibitem{sum of two squares} M. Newmann, \textit{Integral Matrices}, A Series of Monographs and Textbooks, Vol No. $45$, Pure and Applied Mathematics, $(1972)$.

\bibitem{encrypt2}M. Zeriouh, A. Chillali and A. Boua, \textit{Cryptography Based on the Matrices}, Bol. Soc. Paran. Mat., Vol $37$, pp $75-83$, $(2019)$.

\bibitem{} MRS. F. J. MacWilliams, \textit{Orthogonal Circulant Matrices over Finite Fields, and How to Find Them}, Journal of Combinatorial Theory, No. $ 10$, pp $1-17$, $(1971)$.

	\bibitem{}P. Delsarte, J. M. Goethals, and J. J. Seidel, \textit {Orthogonal Matrices with  Zero Diagonal- II}, Canadian Journal of Mathematics, No. $5$, pp $816-832$, $(1971)$. 
	
\bibitem{encrypt1}	R. Harkins, E. Weber, and A. Westmeyer, \textit{Encryption Schemes using Finite Frames and Hadamard arrays}, Experimental Mathematics, Vol. $14$, No. $4$, $(2005)$.
\bibitem{Farrell}P. G. Farrell, \textit{ Linear Binary Anticodes}, Electronics Letters No. $6$, pp $419-421$, $(1970)$

	\bibitem {Pal}R. E. A. C. Paley, \textit{On orthogonal matrices}, Journal of Mathematical Physics, No. $12$, pp $311-320$, $(1933)$.
	
	\bibitem{} Sophie Huczynska (with changes by Max Neunhoffer), \textit{Finite Fields}, Semester $2$, Academic Year $(2012/13)$.
	
\bibitem{}	William Stein, \textit{Elementary Number Theory:	Primes, Congruences, and Secrets}, January $23, (2017)$

	
\end{thebibliography}
\end{document}